\documentclass[10pt]{article}
\usepackage{amsthm}
\usepackage{ams}
\usepackage[all]{xy}

\theoremstyle{plain}
\newtheorem{thm}{Theorem}[section]
\newtheorem{cor}[thm]{Corollary}
\newtheorem{lem}[thm]{Lemma}

\theoremstyle{definition}
\newtheorem{defn}[thm]{Definition}

\theoremstyle{remark}
\newtheorem{rem}{Remark}[thm]

\title{On reduction maps and support problem in $K$-theory and abelian varieties}
\author{Stefan Bara\'{n}czuk \footnote{\textsc{Department of Mathematics, Adam Mickewicz University, Pozna\'{n}, Poland}; \textit{e-mail address}: stefbar@amu.edu.pl}}
\date{March 1, 2005}

\begin{document}

\maketitle

\section{Introduction.}

By $\mathrm{Supp}(m)$ we will denote the set of prime numbers dividing a positive number $m$. P\'{a}l Erd\"{o}s asked the following question:
\begin{quote} Suppose that for some integers $x, y$, the following condition holds
$$ \mathrm{Supp}(x^{n}-1)=\mathrm{Supp}(y^{n}-1) $$ 
for every natural number $n$. Is then $x=y$?
\end{quote}
C. Corrales-Rodrig\'{a}\~{n}ez and R. Schoof answered the question and proved its analogue for number fields and for elliptic curves in \cite{C-RS}.\\
A. Schinzel proved the support problem for the pair of sets of positive integers in \cite{S}. \\
G. Banaszak, W. Gajda and P. Kraso\'{n} examined the support problem for abelian varieties for which the images of the $l$-adic representation is well controled and for K-theory of number fields in \cite{BGK1} and \cite{BGK2}.\\
The support problem for abelian varieties over number fields was considered independently by Ch. Khare and D. Prasad in \cite{KP}.\\
M. Larsen in \cite{Larsen} gave a solution of the support problem for all abelian varieties over number fields.\\

In this paper we consider the generalization of the support problem for $K$-theory and abelian varieties; namely we deal with the pair of sets of points instead of pair of points.
The main technical result of the paper is Theorem \ref{maintheorem} and its corollaries: Theorems \ref{maintheorem K-theory} and  \ref{maintheorem AbVar} which let us control the images of lineary independent points of $K$-groups and abelian varieties over number fields via reduction maps. These theorems are the refinement of Theorem 3.1 of \cite{BGK3} and are proven using similar methods. I have recently found out that R. Pink has proven by a different method a result similar to Theorem \ref{maintheorem AbVar} cf. \cite{Pink}, Corollary 4.3.\\

\textit{Acknowledgments}. This paper forms a part of my PhD thesis. I wish to thank my thesis advisor, Grzegorz Banaszak, for suggesting the problem.\\
I would like to thank Wojciech Gajda for discussions and especially for pointing out an argument in the proof of Lemma 5 in \cite{KP}.\\
The research was partially financed by AAG EU grant MRTN-CT 2003-504917.

\section{Groups of the Mordell-Weil type.}

The following axiomatic setup of Mordell-Weil systems was developed in \cite{BGK3}.\\

\textit{Notation.}\\

 \begin{tabular}{ll}
 $l$		& a prime number \\
 $F$		& a number field, $\mathcal{O}_{F}$ its ring of integers\\
 $G_{F}$	& = $G(\bar{F}/F)$\\
 $v$		& a finite prime of $\mathcal{O}_{F}$ , $\kappa _{v}$ =
$\mathcal{O}_{F} / v$ the residue field at $v$\\
 $g_{v}$	& =  $G(\bar{\kappa _{v}}/\kappa _{v})$ \\
 $T_{l}$	& a free $\mathbb{Z}_{l}$ -module of finite rank $d$\\
 $V_{l}$	& = $T_{l} \otimes _{ \mathbb{Z}_{l} } \mathbb{Q}_{l}$\\
 $A_{l}$	& = $V_{l} / T_{l}$\\
 $S_{l}$	& a fixed finite set of primes of  $\mathcal{O}_{F}$
containing all primes above $l$\\
 $\rho _{l} : G_{F} \to GL(T_{l})$ & a Galois representation unramified
outside the set $S_{l}$ \\
 $\bar{\rho _{l^{k}}} : G_{F}
\to GL(T_{l}/l^{k}) $	& the residual representation induced by $\rho _{l}$\\
 $F_{l^{k}}$	& the number field $\bar{F}^{\textrm{Ker}\bar{\rho
_{l^{k}}}}$\\
 $F_{l^{\infty}}$	& = $\bigcup_{k}^{}{F_{l^{k}}} $\\
 $G_{l}$	&  = $ G(F_{l}/F) $\\
 $G_{l^{k}}$	&  = $ G(F_{l^{k}}/F) $\\
 $G_{l^{\infty}}$	&  = $ G(F_{l^{\infty}}/F) $\\
 $C[l^k]$	&the subgroup of $l^k$ - torsion elements of an abelian
group $C$\\
 $C_{l}$	& =   $\bigcup_{k}^{}{C[{l^{k}}]} $ , the $l$-torsion
subgroup of $C$
 \end{tabular}
 \\

Let $L / F$ be a finite extension and $w$ a finite prime in $L$. We write $w \notin S_{l}$ to indicate that $w$ is not over any prime in $S_{l}$.\\
Let $\mathcal{O}$ be a ring with unity, free as $Z$-module, which acts on $T_{l}$ in such a way that the action commutes with the $G_{F}$ action. All modules over the ring $\mathcal{O}$ considered in this paper are left $\mathcal{O}$-modules.\\
Let $\left\{B(L)\right\}_{L}$ be a direct system of finitely generated $\mathcal{O}$-modules indexed by all finite field extensions $L/F$. We assume that for every embedding $L \hookrightarrow L'$ of extensions of $F$, the induced structure map $B(L) \to B(L')$ is a homomorphism of $\mathcal{O}$-modules.\\
Similarly, for every prime $v$ of $F$ we define a direct system $\left\{B_{v}(\kappa_{w})\right\}_{\kappa_{w}}$ of $\mathcal{O}$ -modules where $\kappa_{w}$ is a residue field for a prime $w$ over $v$ in a finite extension $L/F$ .


We have the following assumptions on the action of the $G_{F}$ and
$\mathcal{O}$:

\begin{description}

\item{(A{1})} for each $l$, each finite extension $L/F$ and any prime $w$
of $L$, such that $w\notin S_{l}$ we have $T_{l}^{Fr_{w}} = 0$ , where
$Fr_{w} \in g_{w}$ denotes the arithmetic Frobenius at $w$;
\item{(A{2})} for every $L$ and $w\notin S_{l}$ there are natural maps
$\psi _{l, L}$, $\psi _{l, w}$ and $r_{w}$ respecting $G_{F}$ and
$\mathcal{O}$ actions such that the diagram commutes:

 \begin{equation}
 \xymatrix{ B(L) \otimes \mathbb{Z}_{l} \ar[r]^{r_{w}} \ar@{^{(}->}[d]_{\psi _{l, L}} &
B_{v}(\kappa _{w})_{l} \ar[d]_{\cong } ^{\psi _{l, w}} \\
 H^{1}_{f, S_{l}}(G_{L}, T_{l}) \ar[r]^{r_{w}} & H^{1}(g_{w}, T_{l})  }\label{axiom A2}
 \end{equation}
where $ H^{1}_{f, S_{l}}(G_{L}, T_{l})$ is the group defined by Bloch and
Kato (\cite{BlochKato}). The left (resp., the right) vertical arrow in the diagram (\ref{axiom
A2}) is an embedding (resp., an isomorphism) for every $L$ (resp., for every
$w \notin S_{l}$);
\item{(A3)} for every $L$ the map $\psi _{l, L}$  is an isomorphism for
almost all $l$ or $B(\bar{F})$ is a discrete $G_{F}$-module divisible by
$l$, for almost $l$;
\item{ (A4)} for every $L$ we have: $B(\bar{F})^{G_{L}} \cong
B(L)$ and $H^{0}(G_{L}, A_{l}) \subset B(L)$.

\end{description}

As in \cite{Ribet}
we impose the following four axioms on the representations which we consider:

\begin{description}

\item{(B1)} $End_{G_{l}}(A_{l}[l])\cong \mathcal{O} / l\mathcal{O}$, for
almost all $l$ and $End_{G_{l^\infty}}(T_{l})\cong \mathcal{O}\otimes
\mathbb{Z}_{l}$, for all $l$;

\item{(B2)} $A_{l}[l]$ is a semisimple $\mathbb{F}_{l}[G_{l}]$-module, for
almost all $l$ and $V_{l}$ is a semisimple
$\mathbb{Q}_{l}[G_{l^\infty}]$-module, for all $l$;

\item{(B3)} $H^{1}(G_{l}, A_{l}[l])=0$, for almost all $l$ and
$H^{1}(G_{l^\infty}, T_{l})$ is a finite group, for all $l$ ;

\item{(B4)} for each finitely generated subgroup $\Gamma \subset B(F)$ the
group
$$\Gamma '= \left\{P \in B(F) : mP \in \Gamma \; \textrm{for some}
\; m \in \mathbb{N}\right\}$$
is such that $\Gamma ' /\Gamma$ has a finite exponent.

\end{description}

For a point $R \in B(L)$ (resp., a subgroup  $\Gamma \subset B(F)$ ) we
denote $\hat{R} = \psi _{l, L} (R)$ (resp., $\hat{\Gamma} = \psi _{l, L}
(\Gamma)$).\\

\begin{defn} The system of modules $\left\{B(L)\right\}_{L}$ fulfilling the above axioms is called the \textbf{Mordell-Weil system} and the groups $B(L)$ are called \textbf{groups of Mordell-Weil type}.
\end{defn}

\section{Examples of Mordell-Weil systems.}\label{Examples of Mordell-Weil systems.}
In both cases below the axioms ($A_{1}$)-($A_{4}$) and ($B_{1}$)-($B_{4}$) are satisfied by \cite{BGK3}, proofs of Theorem 4.1 and Theorem 4.2 loc. cit. In particular, in the case of abelian varieties, assumptions are fulfilled due to results of Faltings \cite{Fa}, Zarhin \cite{Za}, Serre \cite{Serre} and Mordell-Weil.\\

\textbf{Algebraic $K$-theory of number fields.}\\
For every finite extension $L/F$ consider the Dwyer-Friedlander maps (\cite{DwyerFriedlander}):
$$K_{2n+1} \to K_{2n+1}(L)\otimes \mathbb{Z}_{l} \to H^{1}(G_{L} ; \mathbb{Z}_{l}(n+1)) ,$$
where the action of $G_{L}$ on $\mathbb{Z}_{l}(n+1)$ is given by 
the $(n+1)^{\mathrm{th}}$ tensor power of the cyclotomic character.\\
Let $C_{L}$ be the subgroup of $K_{2n+1}(L)$ generated by the $l$-parts of kernels of the maps  $K_{2n+1}(L)  \to H^{1}(G_{L} ; \mathbb{Z}_{l}(n+1)) $ for all primes $l$. By \cite{DwyerFriedlander} $C_{L}$ is finite and according to Quillen-Lichtenbaum conjecture should be trivial.
We put $$B(L)=K_{2n+1}(L)/C_{L}$$ and we have $\mathcal{O} = \mathbb{Z}$. \\

\textbf{Abelian varieties over number fields.}\\
Let $A$ be an abelian variety over number field $F$ and let
$$\rho_{l} : G_{F} \to GL(T_{l}(A))  $$
be the $l$-adic representation given by the action of absolute Galois group on the Tate module $T_{l}(A)$ of $A$. In this case we put $B(L)=A(L)$ for every finite field extension $L/F$ and $\mathcal{O}=End(A)$.


\section{ Kummer theory for $l$-adic representations.}

We introduce the Kummer theory for $l$-adic representations, following \cite{Ribet} and \cite{BGK3}.
Let $\Lambda$ be finitely generated $\mathcal{O}$-submodule of
$B(F)$ generated by the points $P_{1}, \ldots , P_{r}$, which are linearly
independent over $\mathcal{O}$.\\
For natural numbers $k$ we have the Kummer maps:
$$ \phi _{P_{i}}^{k} : G(\bar{F} / F_{l^k}) \to A_{l}[l^{k}]$$
$$ \phi _{P_{i}}^{k}
(\sigma)=\sigma(\frac{1}{l^{k}}\hat{P_{i}})-\frac{1}{l^{k}}\hat{P_{i}} $$
We define $$ \Phi ^{k} :  G(\bar{F} / F_{l^k}) \to
\bigoplus_{i=1}^{r}{A_{l}[l^{k}]}$$
$$  \Phi ^{k}  = (  \phi _{P_{1}}^{k},\ldots ,   \phi _{P_{r}}^{k} ).$$
We define the field:

$$F_{l^{k}}(\frac{1}{l^{k}}\hat\Lambda ) :=\bar{F}^{\textrm{Ker}\Phi ^{k}}$$

Taking the inverse limit in the following commutative diagram
\begin{displaymath}
\xymatrix{ G(\bar{F} / F_{l^k}) \ar[d] \ar[r]^{ \phi _{P_{i}}^{k}} & A_{l}[l^{k}]
\ar[d]^{\times l}\\
G(\bar{F} / F_{l^{k-1}})  \ar[r]^{ \phi _{P_{i}}^{k-1}} & A_{l}[l^{k-1}] }
\end{displaymath}
we obtain a map: $$ \phi _{P_{i}}^{\infty} : G(\bar{F} / F_{l^\infty}) \to
T_{l}.$$
We define: $$ \Phi ^{\infty} :  G(\bar{F} / F_{l^\infty}) \to
\bigoplus_{i=1}^{r}{T_{l}}$$
$$\Phi = (\phi _{P_{1}}^{\infty},\ldots , \phi _{P_{r}}^{\infty}).$$

\begin{lem}\label{intersection} For $k$ big enough

 $$F_{l^{k}}(\frac{1}{l^{k}} \hat \Lambda)\cap F_{l^{k+1}}=F_{l^{k}}.$$
\end{lem}

\begin{proof} The proof follows the lines of Step 1. of the proof of the Proposition 2.2 in \cite{BGK4} that repeats the argument in the proof of Lemma 5 in \cite{KP}. \\
Consider the following commutative diagram:
\begin{displaymath}
\xymatrix{ G (F_{l^{\infty }} (\frac{1}{l^{\infty }}{\hat \Lambda}
)  / F_{l^{\infty }}) \ar[d] \ar[r]  &
			 T_{l}^{r}/l^{m}T_{l}^{r} \ar^{\cong}[d]  \\
	   G(F_{l^{k+1}}(\frac{1}{l^{k+1}}{\hat
\Lambda}) / F_{l^{k+1 }})
\ar[d] \ar[r] & (A_{l}[l^{k+1}])^{r}/l^{m} (A_{l}[l^{k+1}])^{r} \ar^{\cong}[d] \\
G(F_{l^{k}}(\frac{1}{l^{k}}{\hat
\Lambda}) / F_{l^{k }}) \ar[r]
 & (A_{l}[l^{k}])^{r}/l^{m} (A_{l}[l^{k}])^{r}
}
\end{displaymath}
where the horizontal arrows are the Kummer maps and $m \in \mathbb{N}$ is big enough so that  $ l^{m}T_{l}^{r} \subset \mathrm{Im}(G(F_{l^{\infty }} (\frac{1}{l^{\infty }}{\hat \Lambda}
)  / F_{l^{\infty }}) \to T_{l}^{r})$ . Such $m$ exists by \cite{BGK3} Lemma 2.13.\\
Now we see that for $k$ big enough the images of the maps $G(F_{l^{k}}(\frac{1}{l^{k}}{\hat\Lambda}) / F_{l^{k }}) \to (A_{l}[l^{k}])^{r}/l^{m} (A_{l}[l^{k}])^{r} $ must be all isomorphic. Hence the maps $G(F_{l^{k+1}}(\frac{1}{l^{k+1}}{\hat\Lambda}) / F_{l^{k+1 }})\to G(F_{l^{k}}(\frac{1}{l^{k}}{\hat\Lambda}) / F_{l^{k }})$ are surjectiv.\\
Now the diagram
\begin{displaymath}
\xymatrix @rd@C=5pc@R=3pc{ F_{l^{k+1}}(\frac{1}{l^{k+1}}{\hat\Lambda}) \ar@{-}[r] \ar@{-}[d] & F_{l^{k}}(\frac{1}{l^{k}}{\hat\Lambda})  \ar@{-}[d] \\
 F_{l^{k+1}} \ar@{-}[r] &  F_{l^{k}} }
\end{displaymath}
shows that $$F_{l^{k}}(\frac{1}{l^{k}} \hat \Lambda)\cap F_{l^{k+1}}=F_{l^{k}}.$$

\end{proof}

\section{Main technical result.}\label{results}

\begin{thm}\label{maintheorem}
Assume that 
$\rho(G_{F})$  contains an open subgroup of the group of homotheties 
Let $$P_{1}, \ldots , P_{s} \in B(F)$$ be
the points of infinite order, which are linearly independent over
$\mathcal{O}$.\\
Then for any prime $l$, and for any set $\left\{k_{1},\ldots
,k_{s}\right\}\subset \mathbb{N} \cup  \left\{0\right\}$, there are
infinitely many primes $v$, such that the image of the point $P_{t}$ via
the map $$r_{v} : B(F) \to B_{v}(\kappa _{v})_{l}$$ has order equal
$l^{k_{t}}$ for every $t \in \left\{1, \ldots , s \right\}$.

\end{thm}


\begin{proof} 
Let us rename the points $P_{1}, \ldots , P_{s} \in B(F)$ in the following way:
$$P_{1}, \ldots , P_{i}, Q_{1}, \ldots , Q_{j} \in B(F) ,$$ 
and we are going to show that for any prime $l$, and for any set $\left\{k_{1},\ldots
,k_{i}\right\}\subset \mathbb{N} \setminus  \left\{0\right\}$, there are
infinitely many primes $v$, such that the image of the point $P_{t}$ via
the map $$r_{v} : B(F) \to B_{v}(\kappa _{v})_{l}$$ has order equal
$l^{k_{t}}$ for every $t \in \left\{1, \ldots , i \right\}$ and the images of the points $Q_{1}, \ldots , Q_{j}$ are trivial. It is enough to prove the theorem in a case when $k_{t}=1$
for every $t \in \left\{1, \ldots , i \right\}$ .

We will make use of the following diagram:\\

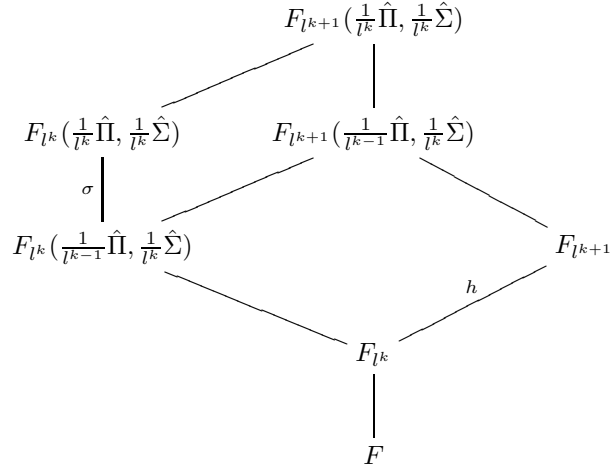
$$\xymatrix{ &F_{l^{k+1}}(\frac{1}{l^{k}}{\hat \Pi} ,\frac{1}{l^{k}}{\hat
\Sigma})\ar@{-}[ld]\ar@{-}[d]&\\
F_{l^{k}}( \frac{1}{l^{k}}{\hat \Pi} ,\frac{1}{l^{k}}{\hat
\Sigma})\ar@{-}[d]_{\sigma}
&F_{l^{k+1}}(\frac{1}{l^{k-1}}{\hat \Pi} ,\frac{1}{l^{k}}{\hat \Sigma})\ar@{-}[dl]\ar@{-}[rd]&  \\
F_{l^{k}}(\frac{1}{l^{k-1}}{\hat \Pi} ,\frac{1}{l^{k}}{\hat \Sigma})\ar@{-}[rd]&&
F_{l^{k+1}}\ar@{-}[ld]_{h}\\ &F_{l^{k}}\ar@{-}[d]&\\
&F& } $$
where $\Pi$ (resp., $\Sigma$) is the $\mathcal{O}$-submodule of $B(F)$
generated by $P_{1}, \ldots , P_{i}$ (resp., by $Q_{1}, \ldots , Q_{j}$).\\
It follows by Lemma \ref{intersection} that for $k$ big enough\\
\begin{displaymath}
F_{l^{k}}(\frac{1}{l^{k-1}}{\hat \Pi} ,\frac{1}{l^{k}}{\hat \Sigma})\cap F_{l^{k+1}}=F_{l^{k}} .
\end{displaymath}\\


\textbf{Step 1.} Consider the following commutative diagram:

\begin{equation}
\xymatrix{ G (F_{l^{\infty }} (\frac{1}{l^{\infty }}{\hat \Pi}
,\frac{1}{l^{\infty }}{\hat \Sigma})  / F_{l^{\infty }}(\frac{1}{l^{\infty }}{\hat \Sigma})) \ar[d] \ar[r]  &
			 T_{l}^{i} \ar@{->>}[d]  \\
	   G(F_{l^{\infty }}(\frac{1}{l^{k}}{\hat
\Pi} ,\frac{1}{l^{\infty }}{\hat \Sigma}) / F_{l^{\infty }}(\frac{1}{l^{\infty }}{\hat \Sigma}))
\ar[d] \ar@{^{(}->}[r] & (A_{l}[l^{k}])^{i}  \ar^{ l}[d] \\
G(F_{l^{\infty }}(\frac{1}{l^{k-1}}{\hat \Pi}
,\frac{1}{l^{\infty }}{\hat \Sigma}) / F_{l^{\infty }}(\frac{1}{l^{\infty }}{\hat \Sigma})) \ar@{^{(}->}[r]
 & (A_{l}[l^{k-1}])^{i}
}\label{sigmadiagram1}
\end{equation}
The horizontal arrows in the diagram (\ref{sigmadiagram1}) are the Kummer maps. The upper horizontal arrow has finite cokernel by \cite{BGK3} Lemma 2.13., so for $k$ big enough 
the horizontal arrows have cokernels bounded independently of $k$. Hence for $k$ big enough there exists $\sigma \in G(F_{l^{\infty }}(\frac{1}{l^{k}}{\hat
\Pi} ,\frac{1}{l^{\infty }}{\hat \Sigma}) / F_{l^{\infty }}(\frac{1}{l^{k-1}}{\hat
\Pi} ,\frac{1}{l^{\infty }}{\hat \Sigma}))$ such that $\sigma$ maps via the Kummer map
\begin{equation}
G(F_{l^{\infty }}(\frac{1}{l^{k}}{\hat
\Pi} ,\frac{1}{l^{\infty }}{\hat \Sigma}) / F_{l^{\infty }}(\frac{1}{l^{k-1}}{\hat
\Pi} ,\frac{1}{l^{\infty }}{\hat \Sigma})) \longrightarrow (A_{l}[l])^{i}
\label{nontrivialsigmainfty}
\end{equation}
to an element whose all $i$ projections on the direct summands $(A_{l}[l])^{i}$ are nontrivial.\\
Then the following tower of fields

\begin{equation}
\xymatrix @rd@C=5pc@R=3pc{ F_{l^{\infty }}(\frac{1}{l^{k}}{\hat \Pi} ,\frac{1}{l^{\infty
}}{\hat \Sigma}) \ar@{-}[r] \ar@{-}[d] & F_{l^{k}}(\frac{1}{l^{k}}{\hat
\Pi} ,\frac{1}{l^{k}}{\hat \Sigma}) \ar@{-}[d] \\
 F_{l^{\infty }}(\frac{1}{l^{k-1}}{\hat \Pi} ,\frac{1}{l^{\infty
}}{\hat \Sigma}) \ar@{-}[r] &  F_{l^{k}}(\frac{1}{l^{k-1}}{\hat
\Pi} ,\frac{1}{l^{k}}{\hat \Sigma}) }
\end{equation}
shows that there exist $\sigma \in G(F_{l^{k }}(\frac{1}{l^{k}}{\hat
\Pi} ,\frac{1}{l^{k}}{\hat \Sigma}) / F_{l^{k }}(\frac{1}{l^{k-1}}{\hat
\Pi} ,\frac{1}{l^{k }}{\hat \Sigma}))$ such that $\sigma$ maps via the Kummer map

\begin{equation}
G(F_{l^{k }}(\frac{1}{l^{k}}{\hat
\Pi} ,\frac{1}{l^{k}}{\hat \Sigma}) / F_{l^{k }}(\frac{1}{l^{k-1}}{\hat
\Pi} ,\frac{1}{l^{k }}{\hat \Sigma})) \longrightarrow (A_{l}[l])^{i}
\label{nontrivialsigma}
\end{equation}
to an element whose all $i$ projections on the direct summands $(A_{l}[l])^{i}$ are nontrivial.\\


\textbf{Step 2.} Let $k$ be big enough that there is an element $\sigma$
constructed in a previous step.\\
From the assumption on the image $\rho _{l}
(G_{F})$ there exists an automorphism $h \in G(F_{l^{k+1}}/F_{l^{k}})$ such
that the action of $h$ on the module $T_{l}$ is given by a nontrivial
homothety $1+l^{k}u_{0}$, for some $u_{0} \in \mathbb{Z}_{l}^{\times}$.\\
We choose an automorphism
$$\gamma \in G(F_{l^{k+1}}(\frac{1}{l^{k}}{\hat \Pi} ,\frac{1}{l^{k}}{\hat
\Sigma})/ F)$$
such that
$$ \gamma \mid _{F_{l^{k}}(\frac{1}{l^{k}}{\hat \Pi} ,\frac{1}{l^{k}}{\hat
\Sigma})} \, = \, \sigma ,$$
$$ \gamma \mid _{F_{l^{k+1}}} \, = \, h .$$
By the Tchebotarev Density Theorem there exist infinitely many prime ideals
$v$ in $\mathcal{O}_{F}$ such that $\gamma$ is equal to the Frobenius
element for the prime $v$ in the extension
$F_{l^{k+1}}(\frac{1}{l^{k}}{\hat \Pi} ,\frac{1}{l^{k}}{\hat \Sigma})/ F$ . In the remainder of the proof we work with prime ideals $v$ we have just selected.\\


\textbf{Step 3.} Using the same argument as in \cite{BGK3}
, Step 4. , we show that $l r_{v}(P_{1}), \ldots , l r_{v}(P_{i}),
r_{v}(Q_{1}), \ldots ,  r_{v}(Q_{j})$ are trivial in $B_{v}(\kappa
_{v})_{l}$ .\\


\textbf{Step 4.} Using similar argument as in \cite{BGK3}
, Step 5. , we show that $ r_{v}(P_{1}), \ldots ,  r_{v}(P_{i}) $ have
order divisible by $l$ in $B_{v}(\kappa
_{v})_{l}$. Hence elements  $ r_{v}(P_{1}), \ldots ,
r_{v}(P_{i}) $  have orders equal $l$ .

\end{proof}


\section{Corollaries of Theorem \ref{maintheorem}.}\label{cormain} 

The following Theorems are the specialization of Theorem \ref{maintheorem} for examples of groups of Mordell-Weil type described in section \ref{Examples of Mordell-Weil systems.}:

\begin{thm}\label{maintheorem K-theory} Let $C_{F}$ denote the subgroup of $K_{2n+1}(F)$ generated by $l$-parts (for all primes $l$) of kernels of the Dwyer-Friedlander map.\\
Let $P_{1}, \ldots , P_{s}$ be nontorsion elements of  $K_{2n+1}(F)/C_{F}$, which are linearly independent over $\mathbb{Z}$.\\ 
Then for any prime $l$, and for any set $\left\{k_{1},\ldots
,k_{s}\right\}\subset \mathbb{N} \cup  \left\{0\right\}$, there are
infinitely many primes $v$, such that the image of the point $P_{t}$ via
the map $$r_{v} : K_{2n+1}(F)/C_{F} \to K_{2n+1}(\kappa _{v})_{l}$$ has order equal
$l^{k_{t}}$ for every $t \in \left\{1, \ldots , s \right\}$.
\end{thm}

\begin{thm}\label{maintheorem AbVar} Let $A$ be an abelian variety defined over number field $F$.\\
Let $P_{1}, \ldots , P_{s}$ be nontorsion elements of  $A(F)$, which are linearly independent over $End(A)$.\\ 
Then for any prime $l$, and for any set $\left\{k_{1},\ldots
,k_{s}\right\}\subset \mathbb{N} \cup  \left\{0\right\}$, there are
infinitely many primes $v$, such that the image of the point $P_{t}$ via
the map $$r_{v} : A(F)\to A_{v}(\kappa _{v})_{l}$$ has order equal
$l^{k_{t}}$ for every $t \in \left\{1, \ldots , s \right\}$.
\end{thm}


\section{Support problem for $K$-theory and abelian varieties.}

Let us return to the notation from section \ref{results}. Let $\mathcal{O}$ be a commutative integral domain.\\

\begin{thm}\label{support}
Let $P_{1}, \ldots , P_{n}, P_{0}, Q_{1}, \ldots , Q_{n}, Q_{0} \in B(F)$
be the points of infinite order.\\
Assume that for almost every prime $l$
the following condition holds in the group $B_{v}(\kappa _{v})_{l}$:\\

For every set of integers $m_{1}, \ldots , m_{n}$ and for almost
every prime $v$ $$m_{1} r_{v}(P_{1})+\ldots + m_{n} r_{v}(P_{n})=
r_{v}(P_{0}) \, \, \mathrm{implies} \, \, m_{1}
r_{v}(Q_{1})+\ldots + m_{n}r_{v}(Q_{n})= r_{v}(Q_{0}).$$ Then
there exist $\alpha _{i}$, $\beta _{i}\in \mathcal{O} \setminus \{
0 \}$ such that $\alpha _{i} P_{i}+\beta _{i} Q_{i}=0$ in $B(F)$
for every $i \in \left\{0, \ldots n\right\}$.
\end{thm}


\begin{proof} Set $m_{i}=0$ for every $i \in \left\{1, \ldots n\right\}$.
We get
\begin{equation}
r_{v}(P_{0})=0 \, \, \mathrm{implies} \, \, r_{v}(Q_{0})=0
\label{implicationsupp}
\end{equation}
for almost every prime $v$. Assume that $P$ and $Q$ are linearly
independent in $B(F)$ over $\mathcal{O}$. By Theorem \ref{maintheorem} there are
infinitely many primes $v$ such that $r_{v}(P_{0})=0$ and  $
r_{v}(Q_{0})$ has order $l$. This contradicts
(\ref{implicationsupp}). Hence there exist $\alpha _{0}$, $\beta
_{0} \in \mathcal{O} \setminus \{ 0 \}$ such that $\alpha _{0}
P_{0}+\beta _{0} Q_{0}=0$ in $B(F)$.\\ Now fix $m_{1}= \ldots =
m_{j-1}=m_{j+1}=\ldots =m_{n}=0$. Let $m_{j}$ be a natural number
such that $m_{j}P_{j}+P_{0} , (m_{j}+1)P_{j}+P_{0}, (m_{j}+2)P_{j}+P_{0}, 
m_{j}Q_{j}+Q_{0}, (m_{j}+1)Q_{j}+Q_{0}, (m_{j}+2)Q_{j}+Q_{0}$ be nontorsion
points. \\ As above we show that there exist $x_{0}, y_{0},
x_{1}, y_{1}, x_{2}, y_{2} \in \mathcal{O} \setminus \{ 0 \}$ such
that \\

\begin{displaymath}
  \left\{ \begin{array}{l}
x_{0}[m_{j}P_{j}+P_{0}] + y_{0}[m_{j}Q_{j}+Q_{0}] = 0 \\
x_{1}[(m_{j}+1)P_{j}+P_{0}] + y_{1}[(m_{j}+1)Q_{j}+Q_{0}] = 0\\
x_{2}[(m_{j}+2)P_{j}+P_{0}] + y_{2}[(m_{j}+2)Q_{j}+Q_{0}] = 0
  \end{array}\right.
\end{displaymath}
Hence
\begin{equation}
  \left\{ \begin{array}{l}
x_{1}y_{0}P_{j} + y_{1}y_{0}Q_{j} = (x_{0}y_{1} -
x_{1}y_{0})[m_{j}P_{j}+P_{0}]\\ 2[x_{2}y_{0}P_{j} + y_{2}y_{0}Q_{j}]
=(x_{0}y_{2} - x_{2}y_{0})[m_{j}P_{j}+P_{0}]
  \end{array}\right.\label{twoeq}
\end{equation}
If $(x_{0}y_{1} - x_{1}y_{0})=0$ or $(x_{0}y_{2} - x_{2}y_{0})=0$ we are done. So assume that
\begin{equation}
 (x_{0}y_{1} -x_{1}y_{0})(x_{0}y_{2} - x_{2}y_{0})\neq 0.
  \label{assump}
\end{equation}
Hence from (\ref{twoeq}) we get
$$y_{0}[x_{0}(x_{2}y_{1}-x_{1}y_{2})+x_{2}(x_{0}y_{1}-x_{1}y_{0})]P_{j}
=
y_{0}[y_{0}(x_{1}y_{2}-x_{2}y_{1})+y_{2}(x_{1}y_{0}-x_{0}y_{1})]Q_{j}$$
If
$y_{0}[x_{0}(x_{2}y_{1}-x_{1}y_{2})+x_{2}(x_{0}y_{1}-x_{1}y_{0})]\neq0$
or
$y_{0}[y_{0}(x_{1}y_{2}-x_{2}y_{1})+y_{2}(x_{1}y_{0}-x_{0}y_{1})]\neq0$
 we are done. So assume that
\begin{displaymath}
  \left\{ \begin{array}{l}
y_{0}[x_{0}(x_{2}y_{1}-x_{1}y_{2})+x_{2}(x_{0}y_{1}-x_{1}y_{0})]=0\\
y_{0}[y_{0}(x_{1}y_{2}-x_{2}y_{1})+y_{2}(x_{1}y_{0}-x_{0}y_{1})]=0.
  \end{array}\right.
\end{displaymath}
Then
\begin{displaymath}
  \left\{ \begin{array}{l}
x_{0}(x_{2}y_{1}-x_{1}y_{2})+x_{2}(x_{0}y_{1}-x_{1}y_{0})=0\\
y_{0}(x_{1}y_{2}-x_{2}y_{1})+y_{2}(x_{1}y_{0}-x_{0}y_{1})=0.
  \end{array}\right.
\end{displaymath}
Hence we get $$(x_{0}y_{1} - x_{1}y_{0})(x_{0}y_{2} - x_{2}y_{0})=
0$$ that contradicts (\ref{assump}). 

\end{proof}

\begin{thm}\label{Pure support}
Let $P_{1}, \ldots , P_{n}, Q_{1}, \ldots , Q_{n}\in B(F)$ be the
points of infinite order.\\ Assume that for almost every prime $l$
the following condition holds in the group $B_{v}(\kappa
_{v})_{l}$:\\

For every set of positive integers $m_{1}, \ldots , m_{n}$ and for
almost every prime $v$ $$m_{1} r_{v}(P_{1})+\ldots + m_{n} r_{v}(P_{n})=0
\, \, \mathrm{implies} \, \, m_{1} r_{v}(Q_{1})+\ldots + m_{n}r_{v}(Q_{n})=
0. $$ Then there exist $\alpha _{i}$, $\beta _{i}\in \mathcal{O} \setminus
\{ 0 \}$ such that $\alpha _{i} P_{i}+\beta _{i} Q_{i}=0$ in $B(F)$ for
every $i \in \left\{1, \ldots n\right\}$. \end{thm}

\begin{proof}  The proof of the theorem is analogous to the proof of the
Theorem \ref{support}:\\
Let $m_{j}$ be a natural number such that
\begin{eqnarray*}
m_{1}P_{1}+\ldots +m_{j-1}P_{j-1} +m_{j}P_{j} +m_{j+1}P_{j+1}+\ldots
+m_{n}P_{n},\\ m_{1}P_{1}+\ldots +m_{j-1}P_{j-1} +(m_{j}+1)P_{j}
+m_{j+1}P_{j+1}+\ldots +m_{n}P_{n},\\ m_{1}P_{1}+\ldots +m_{j-1}P_{j-1}
+(m_{j}+2)P_{j}+m_{j+1}P_{j+1}+\ldots +m_{n}P_{n},\\ m_{1}Q_{1}+\ldots
+m_{j-1}Q_{j-1}+m_{j}Q_{j} +m_{j+1}Q_{j+1}+\ldots +m_{n}Q_{n},\\
m_{1}Q_{1}+\ldots+m_{j-1}Q_{j-1} +(m_{j}+1)Q_{j} +m_{j+1}Q_{j+1}+\ldots
+m_{n}Q_{n},\\ m_{1}Q_{1}+\ldots +m_{j-1}Q_{j-1} +(m_{j}+2)Q_{j}
+m_{j+1}Q_{j+1}+\ldots +m_{n}Q_{n}
\end{eqnarray*}
be nontorsion points.
There exist $x_{0}, y_{0}, x_{1}, y_{1}, x_{2}, y_{2} \in \mathcal{O}
\setminus \{ 0 \}$ such that \\
\begin{displaymath}
  \left\{ \begin{array}{l}
x_{0}[ m_{1}P_{1}+\ldots
+m_{j-1}P_{j-1} +m_{j}P_{j} +m_{j+1}P_{j+1}+\ldots +m_{n}P_{n}] +\\
+y_{0}[
m_{1}Q_{1}+\ldots +m_{j-1}Q_{j-1} +m_{j}Q_{j} +m_{j+1}Q_{j+1}+\ldots
+m_{n}Q_{n}]=0\\

x_{1}[ m_{1}P_{1}+\ldots
+m_{j-1}P_{j-1} +(m_{j}+1)P_{j} +m_{j+1}P_{j+1}+\ldots
+m_{n}P_{n}]+\\
+y_{1}[m_{1}Q_{1}+\ldots +m_{j-1}Q_{j-1} +(m_{j}+1)Q_{j}
+m_{j+1}Q_{j+1}+\ldots +m_{n}Q_{n}]=0\\

x_{2}[m_{1}P_{1}+\ldots
+m_{j-1}P_{j-1} +(m_{j}+2)P_{j} +m_{j+1}P_{j+1}+\ldots +m_{n}P_{n}
]+\\
+y_{2}[m_{1}Q_{1}+\ldots +m_{j-1}Q_{j-1} +(m_{j}+2)Q_{j}
+m_{j+1}Q_{j+1}+\ldots +m_{n}Q_{n}]=0
 \end{array}\right.
\end{displaymath}
The rest of the proof follows the lines of the proof of the Theorem
\ref{support}.

\end{proof}


\begin{rem} Assume that $\mathcal{O}=\mathcal{O}_{E}$ for
some number field $E$.\\ Assume that there exist $\alpha $, $\beta
 \in \mathcal{O}_{E} \setminus \{ 0 \}$ such that $\alpha
P+\beta  Q=0$ in $B(F)$. Then there exist $z \in \mathbb{Z}
\setminus \{ 0 \}$ such that $z\frac{\beta}{\alpha} \in
\mathcal{O}_{E}$ (see: \cite{Mollin}, p.46). Hence $z P+ z\frac{\beta}{\alpha} Q=0$ in
$B(F)$. We can then replace the expression `` $\alpha _{i}$, $\beta
_{i}\in \mathcal{O} \setminus \{ 0 \}$ '' in Theorem \ref{support}
by `` $\alpha _{i}\in \mathbb{Z} \setminus \{ 0 \}$ ,$\beta
_{i}\in \mathcal{O} \setminus \{ 0 \}$ ''.
\end{rem}


\section{Corollaries of Theorems \ref{support} and \ref {Pure support}.}

As in the section \ref{cormain}, we obtain the specializations of Theorems \ref{support} and \ref {Pure support} for $K$-theory and abelian varieties:\\

\begin{thm}\label{support K-theory}
Let $P_{1}, \ldots , P_{s}, P_{0}, Q_{1}, \ldots , Q_{s}, Q_{0} \in K_{2n+1}(F)/C_{F}$
be the points of infinite order.\\
Assume that for almost every prime $l$
the following condition holds:\\

For every set of integers $m_{1}, \ldots , m_{s}$ and for almost
every prime $v$ $$m_{1} r_{v}(P_{1})+\ldots + m_{s} r_{v}(P_{s})=
r_{v}(P_{0}) \, \, \mathrm{implies} \, \, m_{1}
r_{v}(Q_{1})+\ldots + m_{s}r_{v}(Q_{s})= r_{v}(Q_{0}).$$ Then
there exist $\alpha _{i}$, $\beta _{i}\in \mathbb{Z} \setminus \{
0 \}$ such that $\alpha _{i} P_{i}+\beta _{i} Q_{i}=0$ 
for every $i \in \left\{0, \ldots s\right\}$.
\end{thm} 

\begin{thm}\label{Pure support K-theory}
Let $P_{1}, \ldots , P_{s}, Q_{1}, \ldots , Q_{s}\in K_{2n+1}(F)/C_{F}$ be the
points of infinite order.\\ Assume that for almost every prime $l$
the following condition holds:\\

For every set of positive integers $m_{1}, \ldots , m_{s}$ and for
almost every prime $v$ $$m_{1} r_{v}(P_{1})+\ldots + m_{s} r_{v}(P_{s})=0
\, \, \mathrm{implies} \, \, m_{1} r_{v}(Q_{1})+\ldots + m_{s}r_{v}(Q_{s})=
0. $$ Then there exist $\alpha _{i}$, $\beta _{i}\in \mathbb{Z} \setminus
\{ 0 \}$ such that $\alpha _{i} P_{i}+\beta _{i} Q_{i}=0$ in $B(F)$ for
every $i \in \left\{1, \ldots s\right\}$. 
\end{thm}

\begin{thm}\label{support Abelian Varietes}
Let $A$ be an abelian variety defined over number field $F$ such that $End(A)$ is a commutative integral domain.\\
Let $P_{1}, \ldots , P_{n}, P_{0}, Q_{1}, \ldots , Q_{n}, Q_{0} \in A(F)$
be the points of infinite order.\\
Assume that for almost every prime $l$
the following condition holds:\\

For every set of integers $m_{1}, \ldots , m_{n}$ and for almost
every prime $v$ $$m_{1} r_{v}(P_{1})+\ldots + m_{n} r_{v}(P_{n})=
r_{v}(P_{0}) \, \, \mathrm{implies} \, \, m_{1}
r_{v}(Q_{1})+\ldots + m_{n}r_{v}(Q_{n})= r_{v}(Q_{0}).$$ Then
there exist $\alpha _{i}$, $\beta _{i}\in End(A) \setminus \{
0 \}$ such that $\alpha _{i} P_{i}+\beta _{i} Q_{i}=0$ 
for every $i \in \left\{0, \ldots n\right\}$.
\end{thm} 

\begin{thm}\label{Pure support Abelian Varieties}
Let $A$ be an abelian variety defined over number field $F$ such that $End(A)$ is a commutative integral domain.\\
Let $P_{1}, \ldots , P_{n}, Q_{1}, \ldots , Q_{n}\in A(F)$ be the
points of infinite order.\\ Assume that for almost every prime $l$
the following condition holds:\\

For every set of positive integers $m_{1}, \ldots , m_{n}$ and for
almost every prime $v$ $$m_{1} r_{v}(P_{1})+\ldots + m_{n} r_{v}(P_{n})=0
\, \, \mathrm{implies} \, \, m_{1} r_{v}(Q_{1})+\ldots + m_{n}r_{v}(Q_{n})=
0. $$ Then there exist $\alpha _{i}$, $\beta _{i}\in End(A) \setminus
\{ 0 \}$ such that $\alpha _{i} P_{i}+\beta _{i} Q_{i}=0$ for
every $i \in \left\{1, \ldots n\right\}$. 
\end{thm}


\section{The case $ \mathcal{O} = \mathbb{Z} $.}

Let us return to the notation from section \ref{results}. We consider the special case $\mathcal{O} =\mathbb{Z}$. 

\begin{lem}\label{lemma 1} 
Let $P, Q \in B(F)$ be points of infinite order. Assume that for every prime number $l$ the following condition holds in the group  $B_{v}(\kappa
_{v})_{l}$:
 
 For every positive integer $n$ and for almost every prime $v$:\\
\begin{equation}
n r_{v} (P)=0 \, \, \textrm{implies} \, \, n r_{v} (Q)=0 .\label{lemma equation}
\end{equation}
Then there is an integer $e$ such that $Q=e P$.
\end{lem} 
 
\begin{proof}
By Theorem \ref{Pure support} there are $\alpha ,\beta \in \mathbb{Z} \setminus {0}$ such that $\alpha P = \beta Q$. Let $l^{k}$ be the largest power of prime number $l$ that divides $\beta$, $\beta = b l^{k}$. By (\ref{lemma equation}) we have:
$$\alpha r_{v}(P)=0 \, \, \textrm{implies} \, \, \alpha r_{v}(Q)=0 ,$$
hence
$$\beta r_{v}(Q)=0 \, \, \textrm{implies} \, \, \alpha r_{v}(Q)=0$$
and 
$$b l^{k} r_{v}(Q)=0 \, \, \textrm{implies} \, \, \alpha r_{v} (Q)=0 .$$
But obviously $\alpha r_{v}(Q)=0 \, \, \textrm{implies} \, \, b \alpha r_{v}(Q)=0$.
Hence we get:
\begin{equation}
l^{k} r_{v} (b Q)=0 \, \, \textrm{implies} \, \, \alpha r_{v} (b Q)=0 .\label {lemma equation proof}
\end{equation}

By Theorem \ref{maintheorem} there are infinitely many primes $v$ such that the order of $r_{v} (b Q)$ is $l^{k}$. So by (\ref{lemma equation proof}) we get $\alpha r_{v} (bQ) =0$ and $l^{k}$ divides $\alpha$.\\
Now repeating the argument from the proof of the Theorem 3.12 of \cite{BGK3} we show that $Q= \frac{\alpha}{\beta} P$ with $\frac{\alpha}{\beta} \in \mathbb{Z}$.
\end{proof}  

\begin{lem}\label{lemma 2}

Let $P_{1}, P_{2}, Q_{1}, Q_{2} \in B(F)$ be points of infinite order. Assume that for every prime number $l$ the following condition holds in the group  $B_{v}(\kappa
_{v})_{l}$:\\

 For every set of positive integers $m_{1}, m_{2}$ and for almost every prime $v$:\\
\begin{equation}
m_{1} r_{v}(P_{1})+ m_{2} r_{v}(P_{2})=0
\, \, \mathrm{implies} \, \, m_{1} r_{v}(Q_{1})+ m_{2}r_{v}(Q_{2})=
0.\label{lemma 2 equation}
\end{equation}
Then there is an integer $e$ such that $Q_{1}=e P_{1}$ and $Q_{2}=e P_{2}$.
\end{lem} 

\begin {proof}
By Theorem \ref{Pure support} there are integers $\alpha_{1}, \alpha_{2}, \beta_{1}, \beta_{2}$ such that $\alpha _{1} P_{1} = \beta_{1} Q_{1}$, $\alpha _{2} P_{2} = \beta_{2} Q_{2}$. We can assume that $\alpha_{1}, \alpha_{2} >0$.\\
Now we have to consider two cases.\\

First, assume that $P_{1}$ and $P_{2}$ are linearly independent over $\mathbb{Z}$. Hence $P_{1}$ and $\left|b\right|Q_{2}$ are also linearly independent, where $\beta_{2}=b l^{k}$ and $l^{k}$ is the largest power of prime number $l$ that divides $\beta$.\\
By Theorem \ref{maintheorem} there are infinitely many primes $v$ such that $r_{v} (P_{1})=0$ and $r_{v}(\left|b\right|Q_{2})$ has order $l^{k}$.\\
By (\ref{lemma 2 equation}), for $m_{1}=\left|\beta_{1}\right|$ and $m_{2}=\alpha_{2}$, and by the choise of $v$ we have:
$$\left|\beta_{2}\right| r_{v}(Q_{2})=0 \, \, \textrm{implies} \, \, \alpha_{2} r_{v}(Q_{2})=0$$
$$l^{k} r_{v} (\left|b\right|Q_{2})=0 \, \, \textrm{implies} \, \, \alpha_{2}r_{v}(\left|b\right|Q_{2})=0 .$$
Again by the choice of $v$ 
$$\alpha_{2} r_{v}(\left|b\right|Q_{2})=0 .$$
Hence $l^{k}$ divides $\alpha_{2}$. 
Now we repeat again the argument from the proof of the Theorem 3.12 of \cite{BGK3} showing that $Q_{2}=e
_{2} P_{2}$ for some nonzero integer $e_{2}$ and analogously $Q_{1}=e
_{1} P_{1}$ for some nonzero integer $e_{1}$.\\
Now by (\ref{lemma 2 equation}) 
$$r_{v}(P_{1})+ r_{v}(P_{2})=0
\, \, \mathrm{implies} \, \, r_{v}(Q_{1})+ r_{v}(Q_{2})=
0 .$$
Hence 
\begin{equation}\label{lem2proof}
r_{v}(P_{1})+ r_{v}(P_{2})=0
\, \, \mathrm{implies} \, \, (e_{1}- e_{2})r_{v}(P_{2})=
0 .
\end{equation}
Let now $k$ be arbitrary natural number and $l$ be arbitrary prime number. By Theorem \ref{maintheorem} there are infinitely many primes $v$ such that $r_{v} (P_{1}+P_{2})=0$ and $r_{v}(P_{2})$ has order $l^{k}$. Hence by (\ref{lem2proof}), $l^{k}$ divides $e_{1}-e_{2}$. So  $e_{1}-e_{2} =0$.\\

Now we assume that $P_{1}$ and $P_{2}$ are linearly dependent over $\mathbb{Z}$, i.e. there are numbers $x \in \mathbb{N} \setminus \left\{0\right\}$ and $y \in \mathbb{Z} \setminus \left\{0\right\}$ such that $x P_{1}=y P_{2}$. Hence $\alpha_{2}\beta_{1} x Q_{1}=\alpha_{1}\beta_{2} y Q_{2}$.\\
Put $m_{2}=m_{1}\alpha_{1}\beta_{2} y \, \, \textrm{sgn}(\beta_{2} y)$ in (\ref{lemma 2 equation}):
$$
m_{1} r_{v}(P_{1})+ m_{1}\alpha_{1}\beta_{2} y \, \, \textrm{sgn}(\beta_{2} y) r_{v}(P_{2})=0
\, \, \mathrm{implies} \, \, m_{1} r_{v}(Q_{1})+ m_{1}\alpha_{1}\beta_{2} y \, \, \textrm{sgn}(\beta_{2} y)r_{v}(Q_{2})=
0,
$$
hence
$$
m_{1}[1+\alpha_{1}\beta_{2} x \, \, \textrm{sgn}(\beta_{2} y)] r_{v}(P_{1})=0 \, \, \mathrm{implies} \, \, m_{1}[1+\alpha_{2}\beta_{1} x \, \, \textrm{sgn}(\beta_{2} y)] r_{v}(Q_{1})=0
$$
Putting $P:=[1+\alpha_{1}\beta_{2} x \, \, \textrm{sgn}(\beta_{2} y)]P_{1}$, $Q:=[1+\alpha_{2}\beta_{1} x \, \, \textrm{sgn}(\beta_{2} y)]Q_{1}$ we get 
$$ m_{1} r_{v}(P)=0 \, \, \mathrm{implies} \, \, m_{1} r_{v}(Q)=0$$
By Lemma \ref{lemma 1} there is an integer $s$ such that $Q=s P$. So $$ [1+\alpha_{2}\beta_{1} x \, \, \textrm{sgn}(\beta_{2} y)]Q_{1}=s [1+\alpha_{1}\beta_{2} x \, \, \textrm{sgn}(\beta_{2} y)]P_{1} .$$
Hence
$$ Q_{1}=[s(1+\alpha_{1}\beta_{2} x \, \, \textrm{sgn}(\beta_{2} y))-\alpha_{1}\alpha_{2} x \, \, \textrm{sgn}(\beta_{2} y)]P_{1} .$$
i. e. there is an integer $e_{1}$ such that $Q_{1}=e_{1}P_{1}$.\\
Analogically there is an integer $e_{2}$ such that $Q_{2}=e_{2}P_{2}$.\\

 Now, by (\ref{lemma 2 equation}), for $m_{1}=x$, $m_{2}= l^k -y$ where  $k$ is an arbitrary natural number and $l$ is an arbitrary prime number such that $ l^k -y > 0$, we get:
$$ l^k r_{v}(P_{2})=0\, \, \mathrm{implies} \, \, y(e_{1}-e_{2})r_{v}(P_{2})=0 .$$ 
By Theorem \ref{maintheorem} there are infinitely many primes $v$ such that $r_{v} (P_{2})$ has order $l^{k}$, so $l^{k}$ divides $y(e_{1}-e_{2})$. But $k$ was arbitrary, so $e_{1}-e_{2}=0$.
\end{proof}

\begin{thm}\label{Pure support O=Z}
Let $P_{1}, \ldots , P_{n}, Q_{1}, \ldots , Q_{n}\in B(F)$ be the
points of infinite order.\\ Assume that for every prime number $l$
the following condition holds in the group $B_{v}(\kappa
_{v})_{l}$:\\

For every set of positive integers $m_{1}, \ldots , m_{n}$ and for
almost every prime $v$ $$m_{1} r_{v}(P_{1})+\ldots + m_{n} r_{v}(P_{n})=0
\, \, \mathrm{implies} \, \, m_{1} r_{v}(Q_{1})+\ldots + m_{n}r_{v}(Q_{n})=
0. $$ Then there exists an integer $e$ such that $Q_{i}=e P_{i}$ in $B(F)$ for
every $i \in \left\{1, \ldots n\right\}$. \end{thm}

\begin{proof}

There is $s \in \mathbb{N} \setminus \left\{0\right\}$ such that $P:=sP_{2}+P_{3}+\ldots +P_{n}$, $\bar{P}:=(s+1)P_{2}+P_{3}+\ldots +P_{n}$, $Q:=sP_{2}+P_{3}+\ldots +P_{n}$,  $\bar{Q}:=(s+1)Q_{2}+Q_{3}+\ldots +Q_{n}$ are nontorsion points.\\
By the assumption of the theorem the following condition holds for every set of positive integers $m_{1}, m_{2}$ and for
almost every prime $v$:
$$m_{1} r_{v}(P_{1})+ m_{2} r_{v}(P)=0
\, \, \mathrm{implies} \, \, m_{1} r_{v}(Q_{1})+ m_{2}r_{v}(Q)=
0 ,$$
$$m_{1} r_{v}(P_{1})+ m_{2} r_{v}(\bar{P})=0
\, \, \mathrm{implies} \, \, m_{1} r_{v}(Q_{1})+ m_{2}r_{v}(\bar{Q})=
0 .$$
By Lemma \ref{lemma 2} there is an integer $e$ such that $Q_{1}=e P_{1}$, $Q=e P$, $\bar{Q}=e \bar{P}$ , i.e.
 \begin{displaymath}
  \left\{ \begin{array}{l}
e[ sP_{2}+P_{3}+\ldots +P_{n}]=sQ_{2}+Q_{3}+\ldots +Q_{n}\\
e[ (s+1)P_{2}+P_{3}+\ldots +P_{n}]=(s+1)Q_{2}+Q_{3}+\ldots +Q_{n} ,\\
  \end{array}\right. 
\end{displaymath}
hence $Q_{2}=eP_{2}$. Analogically $Q_{i}=e P_{i}$ for
every $i \in \left\{3, \ldots n\right\}$.

\end{proof}

\section{Corollaries of Theorem \ref{Pure support O=Z}.} 

\begin{thm}\label{Pure support K-theory O=Z}
Let $P_{1}, \ldots , P_{s}, Q_{1}, \ldots , Q_{s}\in K_{2n+1}(F)/C_{F}$ be the
points of infinite order.\\ Assume that for every prime $l$
the following condition holds:\\

For every set of positive integers $m_{1}, \ldots , m_{s}$ and for
almost every prime $v$ $$m_{1} r_{v}(P_{1})+\ldots + m_{s} r_{v}(P_{s})=0
\, \, \mathrm{implies} \, \, m_{1} r_{v}(Q_{1})+\ldots + m_{s}r_{v}(Q_{s})=
0. $$ Then there exists $e \in \mathbb{Z} \setminus
\{ 0 \}$ such that $Q_{i}= e P_{i}$ for
every $i \in \left\{1, \ldots s\right\}$. 
\end{thm}

\begin{cor}\label{Schinzel}
Let $p_{1}, \ldots , p_{s}, q_{1}, \ldots , q_{s}\in F^{\ast}$ and suppose that for almost every prime ideal $ \wp$ in $\mathcal{O}_{F}$ and for every set of positive integers $m_{1}, \ldots , m_{s}$ the following condition holds:
\begin{displaymath} \prod_{i=1}^{s} {p_{i}^{m_{i}}}=1 \, \, (\mathrm{mod} \wp) \, \, \, \, \mathrm{implies} \, \, \, \, \prod_{i=1}^{s} {q_{i}^{m_{i}}}=1 \, \, (\mathrm{mod} \wp) .
\end{displaymath}
Then there exists $e \in \mathbb{Z} \setminus
\{ 0 \}$ such that $q_{i}= p_{i} ^{e}$ for
every $i \in \left\{1, \ldots s\right\}$. 
\end{cor}

\begin{rem} A. Schinzel (\cite{S}, Theorem 1.) proved by a different method a similar result. Corollary \ref{Schinzel} is a bit more general since it assumes only positive coefficients $m_{i}$.  
\end{rem} 

\begin{thm}\label{Pure support Abelian Varieties O=Z}
Let $A$ be an abelian variety defined over number field $F$ such that $\mathrm{End}(A)= \mathbb{Z}$.\\
Let $P_{1}, \ldots , P_{n}, Q_{1}, \ldots , Q_{n}\in A(F)$ be the
points of infinite order.\\ Assume that for every prime $l$
the following condition holds:\\

For every set of positive integers $m_{1}, \ldots , m_{n}$ and for
almost every prime $v$ $$m_{1} r_{v}(P_{1})+\ldots + m_{n} r_{v}(P_{n})=0
\, \, \mathrm{implies} \, \, m_{1} r_{v}(Q_{1})+\ldots + m_{n}r_{v}(Q_{n})=
0. $$ Then there exists $e \in \mathbb{Z} \setminus
\{ 0 \}$ such that $Q_{i}= e P_{i}$ for
every $i \in \left\{1, \ldots n\right\}$. 
\end{thm}

\bibliographystyle{plain}

\end{document}